# A *M/G/∞* Queue Busy Period Distribution Characterizing Parameter

Manuel Alberto M. Ferreira[1]

[1] Instituto Universitário de Lisboa, BRU - IUL, Lisboa, Portugal

Correspondence: Manuel Alberto M. Ferreira, Instituto Universitário de Lisboa, BRU - IUL, Lisboa, Portugal. E-mail: manuel.ferreira@iscte.pt

**Abstract**

In this work it is proposed a parameter, $\eta$, to characterize the *M/G/∞* Queue Busy Period distribution that is a modification of the *peakedness* proposed in Whitt (1984). As this distribution moments centered at the origin are given, in general, by very complex expressions it is relevant to search for parameters simpler but with good descriptive properties.

**1. Introduction**

In the *M/G/∞* queue system the customers arrive in accordance with a Poisson process at rate $\lambda$, receive a service which length is a positive random variable with distribution function $G(.)$ and mean $\alpha$ and, upon its arrival, find an available server. Each customer service is independent from the other customers' services and from the arrivals process. The traffic intensity is given by $\rho=\lambda\alpha$.

When operating, in the *M/G/∞* system there is a sequence of busy and idle periods, as in any queue operation. As in the *M/G/∞* busy period distribution the moments centered at the origin are given, with a few exceptions, by very complex expressions it is relevant to search for parameters analytically simpler but with good descriptive properties. The case of this queue with exponentially distributed service times is an example of moments complex expressions strange, because not even the evident markovian properties that it fulfills bring the expected simplifications, and important, because of its use in the applications.

In the sequence of this work it is developed a parameter, called $\eta$, to characterize the *M/G/∞* queue busy period distribution that is a modification of the *peakedness* proposed in the work of Whitt (1984). Its expression, for several service time distributions, in particular the exponential, is quite simple.

This work is organized as follows. In section 2 it is made a state of play on the *M/G/∞* busy period distribution moments. Then, in section 3, it is introduced the *peakedness* for the M/G/∞ busy period distribution and are studied its properties. The parameter $\eta$, the goal of this work, is introduced in section 4, as a result of the *peakedness* recognized properties in the former section. The parameter $\eta$ is thought as an improvement of the *peakedness*. The paper ends with the conclusions and a short list of references.

**2. An Overview on the M/G/∞ Busy Period Moments**

The length of the *M/G/∞* queue busy period is a random variable, designated *B*, with Laplace transform given by, see Stadje (1985):

$$\bar{B}(s) = 1 + \frac{s}{\lambda} - \frac{1}{\lambda \int_0^\infty e^{-st-\lambda \int_0^t [1-G(v)]dv} dt} \qquad (1)$$

Its inversion is a complicate task except in the cases $G_1(t) = \frac{e^{-\rho}}{e^{-\rho}+(1-e^{-\rho})e^{\lambda t}}$, $t \geq 0$, for which *B* is exponentially distributed with a probability concentration, $e^{-\rho}$, at the origin and $G_2(t) = 1 - \frac{1}{1-e^{-\rho}+e^{-\rho+\lambda\frac{1}{1-e^{-\rho}}t}}$, $t \geq 0$, for which *B* is exponentially distributed (Note 1).

For the *M/D/∞* system it is possible to invert $\bar{B}(s)$ (Ferreira, 1996a) but the resulting expression is extremely complex, although allows understanding the probabilistic structure of *B* (Note 2).

As for the moments of *B*, after (1) it is obtained, see (Ferreira & Ramalhoto, 1994):

$$E[B^n] = (-1)^{n+1} \left\{ \frac{e^\rho}{\lambda} nC^{(n-1)}(0) - e^\rho \sum_{p=1}^{n-1} (-1)^{n-p} \binom{n}{p} E[B^{n-p}]C^{(p)}(0) \right\}$$

where,

$$C^{(n)}(0) = \int_0^\infty (-t)^n e^{-\lambda \int_0^t [1-G(v)]dv} \lambda(1-G(t))dt, \ n = 1,2,... \quad (2)$$

So it is possible, with this recurrent formula, to write exact expressions for the moments centered at the origin, obtaining

$$E[B] = \frac{e^\rho - 1}{\lambda}$$

whichever is the service time distribution and

$$E[B^2] = \frac{2e^\rho \int_0^\infty \left(e^{\lambda \int_t^\infty [1-G(v)]dv} - 1\right)dt}{\lambda^2}$$

that depends on the whole service time distribution probabilistic structure (Note 3).

The Expression (2) has evidently no interest in the cases of the $M/G_1/\infty$ and $M/G_2/\infty$ systems. But for the $M/D/\infty$ system

$$C^{(0)}(0) = 1 - e^{-\rho} \text{ and } C^{(n)}(0) = e^{-\rho}(-\alpha)^n - \frac{n}{\lambda}C^{(n-1)}(0), n = 1,2,...$$

and then it is possible to write, by a process of recurrence, simple expressions for the $E[B^n], n = 1,2,...$

So, the situation is easy in the cases of the $M/G_1/\infty$ and $M/G_2/\infty$ systems. For the $M/D/\infty$ system is good only for the moments' computation. In other systems, in particular for the $M/M/\infty$ system, the situation is complicated (Note 4).

In fact it is strange that for the $M/M/\infty$ system, very important in the applications, see for instance (Ferreira, Andrade, Filipe & Coelho, 2011), despite its evident markovian properties, since both the inter-arrival times and the service time are exponentially distributed, it is not possible to have simple exact expressions for the busy period moments.

**3. The *Peakedness***

In Whitt (1984), where it is conducted a search for queues parameters extremal distributions, it is proposed that in the probability distributions related with queues study, the *peakedness*, that is: the value of the Laplace transform in $1/\alpha$, is used to their characterization instead of the variance. That author stands on the works of Holtzman (1973), Eckberg (1977) and Rolski (1972, 1974 and 1976) that make similar proposes for parameters and extremal distributions (Note 5).

For the situation under consideration in this work, designating the *peakedness* by *p*, after (1):

$$p = \bar{B}\left(\frac{1}{\alpha}\right) = 1 + \frac{1}{\rho} - \frac{1}{\lambda \int_0^\infty e^{-\frac{t}{\alpha} - \lambda \int_0^t [1-G(v)]dv} dt} \quad (3)$$

and, also,

$$p = \sum_{n=0}^\infty (-1)^n \frac{E[B^n]}{n!\alpha^n} \quad (4)$$

So, *p* incorporates information about the whole moments of *B*. Then, for the $M/G/\infty$ queue busy period it is poor to consider it a mere option to the variance. It may be considered as a composite parameter that globally characterizes the distribution (Note 6).

For instance:

- **$M/G_1/\infty$** $\qquad\qquad\qquad p^{G_1} = \frac{\rho+1}{e^\rho + \rho}$ $\qquad\qquad$ (5)

- **$M/G_2/\infty$** $\qquad\qquad\qquad p^{G_2} = \frac{\rho}{e^\rho + \rho - 1}$ $\qquad\qquad$ (6)

- **$M/D/\infty$** $\qquad\qquad\qquad p^D = \frac{\rho+1}{e^{\rho+1} + \rho}$ $\qquad\qquad$ (7)

- **$M/M/\infty$** $\qquad\qquad\qquad p^M = \frac{e^\rho - \rho - 1}{\rho(e^\rho - 1)}$ $\qquad\qquad$ (8)

Note that, in these four examples, *p* depends only on $\rho$, being the expressions obtained quite simple. With the exception of $E[B]$ it is the only parameter with a simple expression for the $M/M/\infty$ system *B* distribution, confer with (Ferreira, Andrade, & Filipe, 2012).

For service times distributions related with the exponential, it also holds:

- If the service time distribution is *NBUE* - new better than used in expectation – with mean $\alpha$,

$\int_t^\infty [1 - G(v)]dv \geq \int_t^\infty e^{-\frac{v}{\alpha}}dv$, see (Ross, 1983, p. 273) and, so,

$$p^{NBUE} \leq \frac{e^\rho - \rho - 1}{\rho(e^\rho - 1)} \qquad (9)$$

- If the service time distribution is *NWUE* - new worse than used in expectation – with mean $\alpha$,

$\int_t^\infty [1 - G(v)]dv \leq \int_t^\infty e^{-\frac{v}{\alpha}}dv$, see (Ross, 1983, p. 273) and, so,

$$p^{NWUE} \geq \frac{e^\rho - \rho - 1}{\rho(e^\rho - 1)} \qquad (10)$$

And also, in general,

- If $G(.)$ and $H(.)$ are service time distributions such that $1 - G(t) \leq 1 - H(t)$

$$p^G \geq p^H \qquad (11)$$

## 4. The Parameter $\eta$

Put the Expression (4) in the form

$$p = 1 - \frac{e^\rho - 1}{\rho} + \sum_{n=2}^\infty (-1)^n \frac{E[B^n]}{n!\,\alpha^n}$$

As $1 - \frac{e^\rho - 1}{\rho}$ is insensible to the service time distribution (Note 7) it is defined:

$$\eta = p \frac{\rho}{e^\rho - \rho - 1} + 1 \qquad (12)$$

It was retired from $p$ what does not distinguish the various service time distributions and then the result was normalized dividing for that common part.

After what was seen in the former section, it is possible to obtain simple expressions for $\eta^{G_1}$, $\eta^{G_2}$, $\eta^D$ and $\eta^M$ depending only on ρ. And evidently the properties of $p$ extend to $\eta$.

It is easy to show that

$$1 \leq \eta \leq \frac{e^\rho - 1}{e^\rho - \rho - 1} \qquad (13)$$

for any service time distribution, concluding so that

$$\lim_{\rho \to \infty} \eta = 1 \qquad (14)$$

Then the parameter $\eta$ is computed for some service time distributions and values of $\rho$, see Table 1 below, being *P* the power function with parameter *c=.5* (Note 8).

Table1. Values of $\eta$ for some service times distribution functions

| G | ρ=.5 | ρ=1 | ρ=5 | ρ=10 | ρ=15 | ρ=20 |
|---|---|---|---|---|---|---|
| $G_1$ | 3.3469730 | 1.7488465 | 1.0013731 | 1.0000002 | 1.0000000 | 1.0000000 |
| $G_2$ | 2.4633636 | 1.5121659 | 1.0011518 | 1.0000002 | 1.0000000 | 1.0000000 |
| D | 2.0123054 | 1.3319113 | 1.0005158 | 1.0000001 | 1.0000000 | 1.0000000 |
| M | 2.5414941 | 1.5819767 | 1.0067837 | 1.0000454 | 1.0000000 | 1.0000000 |
| P | 2.4721612 | 1.5747122 | 1.0120525 | 1.0001526 | 1.0000000 | 1.0000000 |

Note that for great values of $\rho$, after $\rho=10$ for some service times distributions and $\rho=20$ for others, under certain conditions, $B$ is approximately exponentially distributed, see (Ferreira & Ramalhoto, 1994) and Ferreira (1996b). In the same conditions $\eta$ assumes the value *1*.

## 5. Conclusions

The parameter $\eta$, defined in Expression (12), incorporates information on the whole moments of $B$ centered at the origin.

It assumes different values when the distributions of $B$ differ as functions of the service distribution.

And there are simple bounds for $\eta$ depending only on the traffic intensity $\rho$.

Even for the *M/M/∞* system, with analytically intractable parameters in the case of the distribution of $B$, $\eta$ is given for a quite simple expression.

It is even admissible that close values of $\eta$, for different service time's distributions, indicate similar behaviors of the respective $B$ distributions.

So $\eta$ is a parameter that can help to characterize the distribution of $B$, discriminating by different service time distributions.

This parameter, $\eta$, may be included in the class of other composite parameters, as for instance the kurtosis or the skewness, that are used to characterize probability distributions putting an emphasis in their probability density functions graphic characteristics. But it is more complete and complex than the two examples here called, since incorporates the whole distribution moments. Because of that, its interpretation in probability density functions graphic characteristics terms is not so obvious and it is an interesting open field of investigation.

**Acknowledgements**

The author thanks very much the work of the Editor and the Reviewers on this paper that resulted in very valuable suggestions for its improvement.

This work was financially supported by FCT through the Strategic Project PEst-OE/EGE/UI0315/2011.

**Notes**

Note 1. Both expressions result from, see (Ferreira & Andrade, 2009),

$$G(t) = 1 - \frac{1}{\lambda} \frac{(1-e^{-\rho})e^{-\lambda t - \int_0^t \beta(u)du}}{\int_0^\infty e^{-\lambda w - \int_0^w \beta(u)du} dw - (1-e^{-\rho})\int_0^t e^{-\lambda w - \int_0^w \beta(u)du} dw}, t \geq 0, -\lambda \leq \frac{\int_0^t \beta(u)du}{t} \leq \frac{\lambda}{e^\rho - 1},$$

making $\beta(t) = \beta$ (constant), being $\beta=0$ for $G_1(t)$ and $\beta = \frac{\lambda}{e^\rho - 1}$ for $G_2(t)$.

Note 2. That is: the *M/G/∞* queue with constant (**D**eterministic) service times.

$$\bar{B}^D(s) = 1 + \frac{1}{\lambda}\left(s - \frac{(s+\lambda)s}{\lambda e^{-(s+\lambda)\alpha} + s}\right)$$

which leads to the probability density function, see Ferreira (1996a)

$$b(t) = \sum_{n=0}^{\infty} \left(\frac{d}{dt}\frac{c(t)}{e^{-\rho}}\right) * \left(\frac{d}{dt}\frac{1-d(t)}{1-e^{-\rho}}\right)^{*n} e^{-\rho}(1-e^{-\rho})^n$$

where * is the convolution operator, $\frac{c(t)}{e^{-\rho}} = \begin{cases} 0, t < \alpha \\ 1, t \geq \alpha \end{cases} = G(t)$ and $\frac{1-d(t)}{1-e^{-\rho}} = \begin{cases} \frac{1-e^{-\lambda t}}{1-e^{-\rho}}, t < \alpha \\ 1, t \geq \alpha \end{cases}$.

Note 3. In Sathe (1985), very simple bounds for the variance of the busy period are presented depending only on $\rho$, $\lambda$ and $\gamma_s$ (service variation coefficient):

$$\lambda^{-2} \max[e^{2\rho} + e^\rho \rho^2 \gamma_s^2 - 2\rho e^\rho - 1; 0] \leq VAR[B] \leq \lambda^{-2}[2e^\rho(\gamma_s^2 + 1)(e^\rho - 1 - \rho) - (e^\rho - 1)^2]$$

Note 4. That is: exponential service times. The *M* to designate the exponential distribution is due to that it fulfills the **M**emoriless property:

$$P(X > t + y | X > t) = P(X > y)$$

Note 5. For Holtzman (1973) and Eckberg (1983) the *peakedness* is the ratio of the variance to the mean of the steady-state number of busy servers in a *GI/M/∞* system associated to a *GI/M/k* system. And Eckberg (1977) showed that this parameter, together with the mean, is a much better second parameter than the variance to characterize the inter-arrival time or service time distribution for the *GI/M/k* system. Then it is recommended in (Whitt, 1984) to use the *peakedness* instead of the variance. As for a renewal process, to know the ratio of the variance to the mean is equivalent to know the Laplace transform evaluated in *1/α*, follows the definition of the *peakedness* in Equation (3).

Note 6. Many approximations for the equilibrium mean queue length in the GI/G/1 queue are based on the first two moments of the interarrival-time and service-time distributions; see Shanthikumar and Buzakott (1980) and Whitt (1982). In Whitt (1984) it is suggested to consider the set-valued function that maps the four moment parameters into the set of possible values of the mean queue length. Here the parameter *p* is composed of infinite moments.

Note 7. In theory of queues, a quantity is said to be insensible to the service time distribution when it depends on it only through its mean.

Note 8. $\alpha = \frac{c}{c+1}$. In this case, $\eta^P$ was computed using directly the Expression (3).